\documentclass[11pt]{amsart}
\usepackage{amsmath,amssymb,amsfonts,amscd,a4}
\usepackage[cmtip, all]{xy}
\usepackage{pstricks}
\usepackage{color}
\usepackage{ulem}  
\usepackage{color} 
\usepackage{verbatim} %%% Added by Shiho 

\theoremstyle{plain}
\newtheorem{thm}{Theorem}[section]

\newtheorem{corollary}[thm]{Corollary}

\theoremstyle{definition}

\numberwithin{equation}{section}
\newcommand{\sD}{{\mathcal D}}
\newcommand{\sE}{{\mathcal E}}
\newcommand{\sF}{{\mathcal F}}

\newcommand{\sO}{{\mathcal O}}

\newcommand{\C}{{\mathbb C}}

\newcommand{\F}{{\mathbb F}}

\newcommand{\Q}{{\mathbb Q}}

\newcommand{\Z}{{\mathbb Z}}

\title [Some fundamental groups in arithmetic geometry]{Some fundamental groups in arithmetic geometry}
\author{H\'el\`ene Esnault} 
\address{Freie Universit\"at Berlin, Arnimallee 3, 14195, Berlin,  Germany}
\email{esnault@math.fu-berlin.de}

\thanks{Supported by  the Einstein program. }

\date{\today} %% changed by Shiho  

\begin{document}

\begin{abstract}  
We report on Deligne's finiteness theorem for $\ell$-adic representations on smooth varieties defined over a finite field, on its crystalline version, and on how  the  geometric \'etale fundamental group of a smooth projective variety defined over a field of positive  characteristic  controls crystals on the infinitesimal site and should control those on the crystalline site.
\end{abstract}
\maketitle

\section{Acknowledgements}\label{ack}
We thank the organizers of the 2015 Summer research Institute on Algebraic Geometry for the invitation to give a series of three plenary lectures on  ``Some Fundamental Groups in Arithmetic Geometry''. 
\medskip

My understanding of the mathematics presented here has been influenced in an important way by discussions and joint work  on the topic with many mathematicians. I thank  them all, notably 
T. Abe, Y. Andr\'e, J.-B. Bost, B. Bhatt, J. de Jong, P. Deligne, V. Drinfeld,  E. Hrushovski, M. Kerz, L. Kindler, A. Langer, {{\textdagger} V. Mehta}, T. Saito, A. Shiho, V. Srinivas, Y. Varshavsky. I also thank the two referees for a thorough reading and for remarks.

\section{ Deligne's conjectures: $\ell$-adic theory}
The classical Hermite-Minkowski theorem asserts that there are finitely many numbers fields with bounded discriminant.  In this section we present Deligne's program aiming at showing an analog
finiteness theorem on complex varieties for variations of Hodge structures and on varieties over finite fields for $\ell$-adic sheaves.

\begin{thm}[Deligne, \cite{Del84}, Thm. 0.5.] \label{thm:finitenessC}
Let $X$ be a complex smooth connected variety, let $ r, w$ be natural numbers with $r\neq 0$. Then there are finitely many rank $r$ $\Q$-local systems which are direct factors of a $\Q$-variation of polarizable pure Hodge structure of  weight $w$.

\end{thm}
The inspiration for this finiteness theorem in Hodge theory comes from   Faltings's finiteness theorem \cite[Cor.~p.344]{Fal83}  for abelian schemes, that is in weight $w=1$.
We refer e.g. to \cite[Section~2]{Cat14} for the notion of a variation of Hodge structure. 
 In particular, as it is regular singular at infinity, the analog of  the discriminant appearing in Hermite-Minkowski's theorem  is just the reduced divisor at infinity in a good normal crossings compactification $X\hookrightarrow \bar X$. This explains why it is enough to fix $X$ and no multiplicities along  the components of $\bar X\setminus X$.

\medskip

On the other hand, for varieties defined over finite fields $\F_q$, in $\ell$-adic theory, $\ell$ prime to $q$, one has the following finiteness theorem.
\begin{thm}[Deligne, \cite{EK12},  Thm.~2.1, \cite{Esn16}, Thm.~3.1]\label{thm:finitenessl} Let $X$ be a normal separated  scheme of finite type defined over a finite field $\F_q$, let  $0\neq r$ be a natural number. Let $D\subset \bar X$ be an effective Cartier divisor of a normal compactification  $\bar X$ with support $\bar X \setminus X$. Then there are finitely many isomorphism classes of irreducible Weil (resp. \'etale) 
rank $r$   lisse $\bar \Q_\ell$-sheaves with ramification bounded by $D$, up to twist with Weil (resp. \'etale) characters of $\F_q$.  The number does not depend on the choice of $\ell$.
\end{thm}
We refer to \cite[Section~1.1]{Del80} for $\ell$-adic sheaves and to \cite[1.1.7]{Del80} for Weil sheaves, also to \cite[Notations]{EK11}.
Here `ramification bounded by $D $' means that the Swan conductor on the pull-back of the sheaf on any smooth curve mapping non-trivially to $X$ is bounded by the pull-back of $D$ (see \cite[Definition~3.6]{EK12}). One could formulate the finiteness theorem by replacing this notion by the one used by Drinfeld in \cite[Thm.~2.5~(ii)]{Dri12}, counting the isomorphism classes of the  sheaves which become tame on a given finite \'etale cover $X'\to X$. That this assumption is stronger is proved in \cite[Proof~of~Prop.~3.9]{EK12}. 

\medskip
Deligne's proof relies  on Lafforgue's main  theorem  which in particular implies  Theorem~\ref{thm:finitenessl} over a smooth curve (\cite[Thm.~VII.6]{Laf02}). The whole question,  how to reduce to curves, is geometric. Some of the key ideas of the proof go back to Wiesend  (\cite{Wie06} and \cite{Wie07}).

\medskip

The chronology is a bit intricate. First  Deligne gave a direct proof of  \ref{cor:numberfield}. This   can be understood as a corollary of Theorem~\ref{thm:finitenessl} which itself was proved later  (see \cite[Thm.~8.2]{EK12} for the deduction).
\begin{corollary}[Deligne, \cite{Del12}, Thm.~3.1,  Deligne's conjecture (ii) in \cite{Del80}, 1.2.10]
\label{cor:numberfield}
Given a  lisse \'etale  $\bar \Q_\ell$-sheaf $V$ with determinant of finite order, the  subfield of $\bar \Q_\ell$ spanned by the coefficients of the minimal polynomials of the Frobenii $F_x$ at closed points $x\in |X|$ acting on $V_{\bar x}$ is a number field. 
\end{corollary}
Using  this in an essential way, Drinfeld proved the existence of  $\ell$-adic companions on smooth quasi-projective varieties defined over a finite field \cite[Thm.~1.1]{Dri12}, which is part (v) of the conjecture \cite[1.2.10]{Del80}.  Then, using Drinfeld's theorem in an essential way, Deligne proved Theorem~\ref{thm:finitenessl} under the additional  assumption that $X$ is smooth.  Finally,
a simple reduction of the problem to the smooth locus of $X$ enables one to extend the finiteness theorem to the case where $X$ is normal (\cite[Thm.~3.1]{Esn16}).

\medskip

Fixing a good compactification $X \hookrightarrow \bar X$, with a strict normal crossings divisor at infinity, then a curve $\bar C$, complete intersection of ample divisors in $\bar X$ in good position,  fulfils the Lefschetz theorem on topological fundamental groups, that is the homomorphism
$$\pi_1^{\rm top}(C:=X\cap \bar C)\to \pi_1^{\rm top}(X)$$
is surjective. In particular,  this reduces Theorem~\ref{thm:finitenessC} to the case where  $X$ is of dimension $1$. 
However, for $X$ of dimension $\ge 2$ in characteristic $p>0$, there is no  Lefschetz theorem.
 Thus Theorem~\ref{thm:finitenessl}  does  not  obviously reduce  to  dimension $1$.
 
 \medskip

 All one has at disposal are two kinds of Lefschetz theorems, one for reducing all  tame coverings of $X$ to  one well chosen curve, one for reducing one specific object ($\ell$-adic representation or given Galois cover) to one curve adapted to this object.
 
  \begin{thm}[Drinfeld, \cite{Dri12}, Prop.~C.2, \cite{EK15}, Section~6] \label{thm:drinfeld}
Let  $\bar X\supset X$ be a projective normal  geometrically connected compactification of a smooth scheme of finite type $X$ defined over a field $k$,  let $\Sigma\subset \bar X$ be a closed subset of codimension $\ge 2$ such that $(\bar X\setminus \Sigma)$ and $(\bar X\setminus \Sigma ) \cap (\bar X \setminus X)$ are smooth. Let $\bar C\subset \bar X\setminus \Sigma$ be a smooth projective curve, complete intersection of ample divisors, meeting $\bar X\setminus X$ transversally. Then  the restriction to $C=\bar C \cap X$ of any finite \'etale connected cover of $X\setminus \Sigma$, which is tame along 
$(\bar X\setminus \Sigma ) \cap (\bar X \setminus X)$,
  is connected. In particular, the homomorphism on the tame fundamental groups
  $\pi^t_1(C)\to \pi_1^t(X)$ is surjective.

\end{thm}
 An important point is that one does not need a good compactification for Drinfeld's theorem. If one has one, one can enhance the theorem to a complete version of  the Lefschetz theorems
 under the Lefschetz conditions    $Lef(X,Y)$ (formal sections along $Y$  of vector bundles   lift to an open neighbourhood of $Y$) and 
 under the effective Lefschetz conditions $Leff(X,Y)$  (formal bundles along $Y$  lift to an open neighbourhood of $Y$)  (see \cite[Thm~2.5]{EK15}, and \cite[X.2,p.90]{SGA2} for Grothendieck's Lefschetz and effective Lefschetz conditions).  
 Using Theorem~\ref{thm:drinfeld},   together with the existence  alterations  \cite{dJ97}, one can prove   Theorem~\ref{thm:finitenessl} with the stronger assumption on the ramification being killed  by one fixed finite \'etale cover $X'\to X$ 
 purely geometrically, without using the existence of $\ell$-adic companions (see \cite[Thm.~1.4]{Esn16}). 
 
 \medskip
 
 For non-tame $\ell$-adic sheaves or covers, only  a much weaker version of the Lefschetz theorems is available.
 \begin{thm}[see e.g. \cite{EK12}, Prop.~B.1, Lem.~B.2] \label{thm:bertini}
Let $X$ be a smooth quasi-projective variety defined over $\F_q$, let $S\subset |X|$ be a finite set of closed points. 
\begin{itemize}
\item[1)] Let $V$  be an irreducible $\bar \Q_\ell$-Weil or -\'etale  lisse sheaf,   then there is a smooth curve $C\to  X$ with $S\subset |C|$, such that $V|_C$ is irreducible.
\item[2)]  Let $H\subset \pi_1(X)$ be an open normal subgroup,  then there is a smooth curve $C\to  X$ with $S\subset |C|$, such that the homomorphism $\pi_1(C)\to \pi_1(X)/H$ is surjective.
\end{itemize}
\end{thm}
 Nonetheless, it has the important following consequences.
 \begin{corollary}[Drinfeld \cite{Dri12}, Thm.~1.1,  Deligne's conjecture (v) \cite{Del80} 1.2.10]
\begin{itemize}
\item[1)] If $V$ is an irreducible Weil sheaf, such that 
 ${\rm det}(V)$ is of finite order, then $V$ has weight $0$.

\item[2)]
If $V$ is an irreducible Weil lisse $\bar \Q_\ell$-sheaf with determinant of finite order, and $\sigma \in {\rm Aut}(\bar \Q_\ell/\Q)$, 
there is an irreducible Weil lisse $\bar \Q_\ell$-sheaf $V_\sigma$, called  the $\sigma$-{\rm companion} of  $V$, with determinant of finite order, such that the characteristic polynomials $f_V, f_{V_\sigma}\in \bar \Q_\ell[t]$ of the local Frobenii $F_x$ satisfy $f_{V_\sigma}=\sigma(f_V)$.
 \end{itemize}
\end{corollary}
In 1) and 2), $V$ and $V_\sigma$ are in fact    \'etale by \cite[Thm.~1.3.1]{Del80}.
Deligne's finiteness Theorem~\ref{thm:finitenessl} for rank $1$ sheaves can be proven directly. 
\begin{thm}[Kerz-Saito, \cite{KS14}, Thm.~1.1] \label{thm:KS}
Let  $X$ be a smooth quasi-projective variety over a perfect field $k$, let $X\subset \bar X$ be a projective compactification with simple normal crossings at infinity,  let $D$ be an effective divisor with support in $ \bar X\setminus X$. Define $\pi_1^{\rm ab}(X, D)$ by the condition that a character $\chi: \pi_1(X)\to \Q/\Z$ factors through 
$\pi_1^{\rm ab}(X, D)$ if and only if the Artin conductor of  $\chi$ pulled-back to any curve $C\to X$ is bounded by the pull-back of $D$ via $\bar C\to \bar X$, where $\bar C$ is a compactification of $C$, smooth along $\bar C\setminus C$. Then the full package of the Lefschetz  theorems holds:  

\medskip 
\noindent 
for a  sufficiently ample  divisor $i: \bar Y\subset \bar X$  in good position with respect to  $\bar X\setminus X$, the homomorphism $$i_*: \pi_1^{\rm ab}(Y, \bar Y\cap D)\to 
\pi_1^{\rm ab}(X, D)$$ is an isomorphism if dim $Y \ge 2$, and is  surjective if dim $Y=1$, where $Y=\bar Y\cap X$. In particular, if  $k=\F_q$, then $${\rm Ker} \big(\pi_1^{\rm ab}(X, D)\to \pi_1^{\rm ab}(k)\big)$$ is finite.
\end{thm}
Theorem~\ref{thm:KS} implies the rank $1$ case of Deligne's finiteness  Theorem~\ref{thm:finitenessl}  in case one has a good compactification. In fact, one does  not need the full package,  only that if $\bar C$ is a complete intersection curve of such hypersurfaces $\bar Y$ as in the theorem, then $$\pi_1^{\rm ab}(C, \bar C \cap D)\to 
\pi_1^{\rm ab}(X, D)$$
is surjective. So far, one does not have tools  to understand a version of this for the whole fundamental group, which  would explain Theorem~\ref{thm:finitenessl} in general.

\section{Deligne's conjectures: crystalline theory}
Let $X$ be a smooth  geometrically connected scheme of finite type over a perfect field $k$ of characteristic $p>0$, $W:=W(k)$ be the ring of Witt vectors, $K={\rm Frac}(W)$ be its field of fractions.  We refer to \cite[Section~1]{ES15} for the following presentation. 

\medskip

One defines the  {\it  crystalline sites} $X/W_n$ as PD-thickenings  $(U\hookrightarrow T/W_n, \delta)$,  where the coverings come from $U\subset X$ Zariski open. The crystalline site 
$X/W$ is then  the $2$-inductive limit of the $X/W_n$, see  \cite[Ch.~7,p.~7-22]{BO78}.
The category of {\it crystals}   ${\rm Crys}(X/W)$  is the category of sheaves of $\sO_{X/W}$-modules of finite presentation, with  transition maps  being  isomorphisms. It is  $W$-linear.  
The category of {\it isocrystals}   ${\rm Crys}(X/W)_\Q$ is its   $\Q$-linearisation.   It is $K$-linear, tannakian.  

\medskip

The absolute Frobenius $F$ acts on  ${\rm Crys}(X/W)_\Q$.  The  largest full subcategory 
${\rm Conv}(X/K) \subset {\rm Crys}(X/W)_{\Q}$ on which every object is 
$F^\infty$-divisible is 
the $K$-tannakian  subcategory of {\it convergent} isocrystals (Berthelot-Ogus)  (Ogus defines the site of enlargements from $X/W$, then convergent isocrystals are crystals  on it of $\sO_{X/K}$-modules of finite presentation).

\medskip

We introduce various categories of $F$-{\it isocrystals}. One defines the category
 $F$-${\rm Conv}(X/K) $  of {\it convergent} $F$-{\it isocrystals} 
as pairs $(\sE, \Phi)$ where $\sE \in 
{\rm Conv}(X/K)$ and 
$ \Phi: F^*\sE\xrightarrow{\cong} \sE$ is a $K$-linear isomorphism. Its is a 
 $\Q_p$-linear tannakian category. The category $F$-${\rm Overconv}(X/K)$ of {\it overconvergent} $F$-{\it isocrystals}, defined analytically by Berthelot in (\cite[2.3.6,~2.3.7]{Ber96} (see also  \cite[p.~288]{LeS07}),
 has a more algebraic description due to Kedlaya. 
 It  consists of those convergent $F$-isocrystals $\sE$ which have unipotent local monodromy after alteration in the sense of Kedlaya (\cite[Introduction]{Ked04}, \cite[Introduction~and~Section~3.2]{Ked07}).  It is a $\Q_p$-linear category, fully embedded in  $F$-${\rm Conv}(X/K) $ (\cite[Thm.~1.1]{Ked04}). We shall just need that if $X$ is proper, then Kedlaya's full embedding is an equivalence, 
 and that
 the group of extensions of the trivial object by itself in this category is $H^1_{\rm rig}(X/K)$, the first rigid cohomology group.
 
 \medskip

When $k$ is a finite field $\F_q, \ q=p^s$, 
one defines $F_{\F_{q}}=F^{s}$-${\rm Overconv}(X/K)$ by the same formulae as before, but now $F^{s}$ acts instead of  $F$. It is a 
 $K$-linear category. Fixing an algebraic closure $\bar{\Q}_p$ of $\Q_p$, one defines 
 $F$-$ {\rm Overconv}(X/K)_{\bar{\Q}_p}$ as its $\bar{\Q}_p$-linearization. It is called  is the category of {\it overconvergent $F$-isocrystals over} $\bar{\Q}_p$. See \cite[1.4.11, 4.1.2]{Abe13}.

 \medskip
 With those notations at hand, one can formulate the general hope.  For smooth geometrically irreducible schemes of finite type $X$ defined over a finite field $\F_q$,
 there should be an analogy   between 
 \begin{itemize}
 \item[i)]
  irreducible objects in $F$-${\rm Overconv}(X/K)_{\bar \Q_p}$ with  determinant of finite order;  
\item[ii)] irreducible lisse
$\bar \Q_\ell$-sheaves  with  determinant of finite order.
\end{itemize}
Upon bounding ramification at infinity in i) and ii), the analogy should extend to 
 irreducible $\Q$-variations of polarisable pure Hodge structures  definable over $\Z$ over  complex varieties. Of course, since there are many categories of isocrystals, one may wonder why those particular ones are the right analogs. 
 The best demonstration is clearly Abe's Theorem~\ref{thm:abe}. But already before the theorem was known, one knew that isocrystals have slopes (a topic not discussed here) and that those isocrystals 
 which come from the variation of crystalline cohomology of the fibers of a smooth projective family have  pure slope parts  which are convergent subisocrystals. However, the lisse $\bar \Q_\ell$-sheaves computing the variation of the $\ell$-adic cohomology of the fibers tend to be irreducible if the geometric variation of the family is big, e.g. if the family is the universal family on a moduli space.

 \begin{thm}[Abe,  \cite{Abe13}, Thm.~4.3.1] \label{thm:abe}
Let $X$ be a smooth curve defined over a finite field  $\F_q$. Then 
\begin{itemize}
\item[1)] an irreducible object in $F$-${\rm Overconv}(X/K)_{\bar \Q_p}$  with  determinant of finite order  is 
$\iota$-pure of weight $0$;
\item[2)]
 an irreducible  lisse  $\bar \Q_\ell$- \'etale sheaf with determinant of finite order has a companion which is an irreducible overconvergent 
 $F$-isocrystal  over $\bar{ \Q}_p$ with  determinant of finite order; vice-versa, an irreducible overconvergent 
 $F$-isocrystal  over $\bar{ \Q}_p$ with determinant of finite order has a companion which is an irreducible  lisse  $\bar \Q_\ell$- \'etale sheaf with  determinant of finite order.
 \end{itemize}
\end{thm}
Here $\iota$ is a fixed isomorphism $\iota: \bar \Q_p\to \C$, and $\iota$-pure means that the Frobenii at all closed points  $x$ of $X$ act on the fiber $E_x $ in 
$$F{\rm-Overconv}(x/{\rm Frac}(W(k(x)))_{\bar \Q_p}=
F{\rm-Conv}(x/{\rm Frac}(W(k(x)))_{\bar \Q_p}
$$
 of the object $\sE$  in $F$-${\rm Overconv}(X/K)_{\bar \Q_p}$ with eigenvalues of complex absolute value $q^{w/2}$ for a fixed real number $w$ called the weight. (Similarly, $\iota$-mixed means that $\sE$ is filtered in $F$-${\rm Overconv}(X/K)_{\bar \Q_p}$ such that the associated graded $gr \sE$ in 
a sum of $\iota$-pure objects. See \cite[Defn.~2.1.3]{AC13}.)

\medskip

Deligne's program \cite[1.2.10]{Del80} in higher dimension on the crystalline side  is not yet achieved. However,  a Lefschetz theorem  such as Theorem~\ref{thm:bertini} for  overconvergent $F$-isocrystals over $\bar{\Q}_p$ has been proven.

\begin{thm}[Abe-Esnault, \cite{AE16}, Thm.~0.1] \label{thm:lefschetzp}
Let $X$ be a smooth connected quasi-projective variety defined over $\F_q$. Let $M$  be an irreducible overconvergent $F$-isocrystal.    
Then there is a open dense subscheme $U\hookrightarrow X$ such that for any finite set $S\subset U$ of closed points, there is a smooth irreducible curve $C\to X$ such that $S\subset C$ and such that $M|_C$ is irreducible.

\end{thm}

 This implies the following.   
\begin{itemize}
\item[1)]  $M$ is  $\iota$-pure of weight $0$ (\cite[Thm.~4.2.2]{Abe13}); 
\item[2)]  Corollary~\ref{cor:numberfield}  remains true with $V$ replaced by $M$, that is 
there is a number field  containing all the coefficients of all local eigenpolynomials (\cite[1.5]{AE16}) at closed points (\cite[Lem.~4.1]{AE16});
\item[3)]    $M$ has $\ell$-adic companions which are irreducible $\bar \Q_\ell$-lisse sheaves (\cite[Thm.~4.3]{AE16});
\item[4)] There is a crystalline version Theorem~\ref{thm:finitenessl}: Let $X$ be a smooth connected quasi-projective variety defined over $\F_q$, $(r, D)$ be as in Theorem~\ref{thm:finitenessl}, $\sigma$ be a field isomorphism from $\bar \Q_p$ to $\bar \Q_\ell$ for some prime number $\ell$ different from  $p$. Then there are finitely many isomorphism classes of 
irreducible overconvergent $F$-isocrystals, up to twist with rank $1$ isocrystals on $\F_q$, such that the $\sigma$-companion (which  by 3) is an  irreducible $\bar \Q_\ell$-lisse sheaf) 
 has ramification bounded by $D$.

\end{itemize}

We note that  in an `unstable preprint'  posted on his webpage, unstable   in the author's terminology, Kedlaya uses weights to deduce 1) and 2) as well as the part of 3) concerning the existence of $\ell$-adic companions. The properties of being lisse and irreducible seems to be inaccessible without the Lefschetz theorem \ref{thm:lefschetzp}. 

\medskip

We also mention \cite[Thm.~1.2]{Kos15} in which a weak analog to Theorem~\ref{thm:finitenessl}
is proven:  if $X$ is a smooth geometrically connected variety defined over a finite field, then a semi-simple unit-root overconvergent $F$-isocrystal in  $F$-${\rm Overconv}(X/K)_{\bar \Q_p}$ is
 isotrivial.  The point is that such an object necessarily is locally isotrivial at infinity, which  reduces the problem  to the case of $X$ smooth projective, thus by the standard Lefschetz  theorem to the curve case. 
  One then applies Abe's theorem \cite[Thm.~4.1]{Abe13} which reduces the statement to Lafforgue's theorem \cite[Thm.~VII.6]{Laf02}).  

\section{Mal\v{c}ev-Grothendieck's theorem, Gieseker's conjecture, de Jong's conjecture}

Let $X$ be a smooth geometrically irreducible scheme of finite type over  field $k$ of characteristic $0$. 
Grothendieck defined the {\it infinitesimal site} $X_\infty$ (\cite{Gro68}) with objects  $U\hookrightarrow T$ where $T$ is an infinitesimal thickening of a Zariski open subscheme $U$, and where coverings  come from the $U$s. {\it Crystals} are finitely presented crystals on $X_\infty$. The category is equivalent to the category of bundles on $X$ with an integrable connection $(E,\nabla)$ or equivalently to the category of $\sO_X$-coherent $\sD_X$-modules. It is a 
$k$-linear category, which is tannakian. 

\begin{thm}[Mal\v{c}ev  \cite{Mal40}, Grothendieck  \cite{Gro70}, Thm.~4.2]  \label{thm:MG}
Let $X$ be a complex smooth variety. If its  \'etale fundamental  group is trivial, then there are no non-trivial  crystals in the infinitesimal site  (with regular singularities at infinity  in case $X$ is not projective).
\end{thm}

Here one uses the Riemann-Hilbert correspondence to translate the assertion on representations of groups which are finitely generated, applied to the topological fundamental group, to the assertion on crystals in the infinitesimal site. The proof then just uses that a $GL(r, A)$-representation, where $A$ is $\Z$-algebra of finite type, is trivial if and only if it is by restriction to the closed points of ${\rm Spec}(A)$.
\medskip

Gieseker \cite[p.~8]{Gie75} conjectured that the analog theorem remains true in characteristic $p>0$. 
Let $X$ be a smooth  projective geometrically irreducible  variety over a field $k$ of characteristic $p>0$. One defines $X_\infty$ and crystals as in characteristic $0$.   By Katz' theorem \cite[Thm.~1.3]{Gie75}, which relies on Cartier descent,  it is also equivalent to the category of Frobenius divisible $\sO_X$-coherent sheaves, that is infinite sequences $(E_0, E_1,\cdots, \sigma_0, \sigma_1, \cdots)$  of bundles $E_n$ on the $n$-th Frobenius twist $X^{(n)}$ of $X$, together with isomorphisms $\sigma_n$ between $E_n$ and the Frobenius pull-back of $E_{n+1}$. Then Gieseker's conjecture predicts that Theorem~\ref{thm:MG} holds in characteristic $p>0$. It has been proved in 2010.

\begin{thm}[Esnault-Mehta, \cite{EM10}, Thm.~1.1] \label{thm:EM}
Let $X$ be a  smooth projective geometrically irreducible variety over a field $k$ of characteristic $p>0$. If its  geometric  \'etale fundamental  group  is trivial, then there are no non-trivial  crystals in the infinitesimal site.
\end{thm}

What in the proof  replaces the finite generation of the topological fundamental group is the existence of quasi-projective moduli for stable bundles with vanishing Chern classes (Langer, \cite[Thm.~4.1]{Lan04}). What then replaces the criterion for triviality is Hrushovski's theorem on the existence of preperiodic points on dominant correspondences over finite fields \cite{Hru04}. Varshavsky in \cite{Var14} gave a proof  of it in the framework of arithmetic geometry, without using  model theory. 

\medskip

One can formulate variants of Gieseker's conjecture.
If  $X$ is not proper, then the  theory of {\it regular singular} crystals in the infinitesimal site has been  developed by Kindler  
\cite{Kin15},
 in such a way that for  those objects with  a finite Tannaka group, it coincides with the notion of {\it tame} quotient of the \'etale fundamental group.
There is  no  good  higher ramification theory  so far, nor does one have an analog of Theorem~\ref{thm:EM},  except for  the tame abelian quotient of the geometric fundamental group  (Kindler, \cite[Thm.~1.4]{Kin13}),  and in the case where $X$ is the smooth locus of a normal projective variety defined over $k=\bar \F_q$   (Esnault-Srinivas \cite[Thm~1.1]{ESB15}; the proof uses Bost's improvement of Grothendieck's LEF theorem,  see \cite[Appendix]{ESB15}).

\medskip

In 2010, de Jong formulated the corresponding conjecture in the category of isocrystals. 
Let $X$ be a smooth  projective geometrically irreducible  variety over a  perfect field $k$ of characteristic $p>0$. If  its  geometric  \'etale fundamental  group  is trivial, then the conjecture predicts that  there are no non-trivial  isocrystals.  As of today, there is no complete understanding of the conjecture. We now list the known results concerning it. 

\medskip
One defines $N(1)=\infty, N(2)=2,N(3)=1, N(r)=1/M(r)$  where for any natural number $r\ge 4$, $M(r)$ is the maximum of the lower common multiples of $a$ and $b$ for all choices of natural numbers $a, b\ge 1$ with $a+b\le r$. For any torsion-free coherent sheaf $\sF$, one denotes by $\mu_{\rm max}(\sF)$ its maximal slope.

\begin{thm}[Esnault-Shiho] \label{thm:ES}
Let $X$ be a smooth  projective geometrically irreducible  variety over a  perfect field $k$ of characteristic $p>0$. 
\begin{itemize}
\item[1)] If the abelian quotient of the geometric \'etale fundamental group of $X$ is trivial, there are no non-trivial isocrystals which are successive  extensions of rank $1$ isocrystals. (\cite[Prop.~2.9, Prop.~2.10]{ES15}). 

\item[2)] If the geometric \'etale fundamental group of $X$ is trivial,  $\mu_{\rm max}(\Omega^1_X)<N(r)$, the irreducible constituents of the Jordan-H\"older filtration of $\sE$ have rank $\le r$, and $\sE$ itself is either  in ${\rm Conv}(X/K)$ or else each of its irreducible constituents has a locally free lattice and has rank $\ge r$, then $\sE $ is trivial. (\cite[Thm.~1.1]{ES15} and \cite[Thm.~1.2]{ES15b}).

\end{itemize}
Let $f: Y\to X$ be a smooth proper  morphism between smooth proper schemes of finite type. 
\begin{itemize}
\item[3)] If the geometric \'etale  fundamental group of $X$ is trivial, then all the Gau{\ss}-Manin convergent isocrystals $R^nf_*\sO_{Y/K}$ are trivial.  If $k$ is finite,  $f$ is projective and $p \geq 3$, one can drop the properness assumption on $X$ (\cite[Thm.~1.3,  Rmk.~1.4]{ES15b}).
\end{itemize}
\end{thm}
We first discuss 3). Le us assume that $k$ is a finite field. Then the statement relies on
\begin{thm}[Abe's \v{C}ebotarev's density theorem, \cite{Abe13}, A.3] \label{thm:cebotarev}
Let $X$ be a smooth scheme of finite type defined over a finite field $k$. If $\sE$ and $\sE'$ are $\iota$-mixed overconvergent $F$-isocrystals over $\bar \Q_p$ with the same set of Frobenius eigenvalues on closed points of $X$, then the semi-simplifications of $\sE$ and $\sE'$ are isomorphic. 

\end{thm}
 
When $X$ is proper or $f$ is projective, 
the  convergent $F$-isocrystal $R^nf_*\sO_{Y/K}$ 
is an overconvergent $F$-isocrystal over $\bar \Q_p$  via the faithful embedding (\cite[Cor.~5.4]{Laz15}), thus obeys Theorem~\ref{thm:cebotarev}.  The Weil conjectures \cite[Thm.~1.1]{KM74}, \cite[Cor.~1.3]{CLS98} in the proper case,
enable one to conclude that the semi-simplification of $R^nf_*\sO_{Y/K}$  is constant. Forgetting the $F$-structure, it is thus a successive extension of the trivial overconvergent isocrystal by itself, thus is trivial, as the first rigid cohomology of $X$ is controlled by the first $\ell$-adic cohomology of $X\otimes_k \bar k$ when $X$ is proper or $p \geq 3$. The latter is trivial   if the geometric fundamental group is trivial. Over a general field $k$, the properness assumption on $X$ allows to compare the statement to the one over finite fields by base change.

\medskip

We discuss 2).  Since the conjecture concerns isocrystals, it is not natural to try to argue with lattices, that is $p$-torsion free crystals  $E$ in a given isocrystal class $\sE$. Unfortunately, there is at present no other way to do, and indeed, basically 2) is proven by showing that under the given assumptions, 
the value  $E_X$ on $X$ of a well chosen crystal  $E$ in the isocrystal class $\sE$ is trivial.  Then one studies the possible liftings modulo $p$-powers. In order to show triviality of $E_X$, one applies Theorem~\ref{thm:EM}. To do so, one has to show the existence of such an $E_X$ which is semi-stable with vanishing numerical Chern classes, so as to be able to define its moduli point. If $E_X$ was $F^\infty$-divisible, then one could apply Theorem~\ref{thm:EM} directly.  This is not the case, even if $\sE\in {\rm Conv}(X/K)$, that is even if $\sE$ is $F^\infty$-divisible. Instead, one shows that $F^a$-divisibility is enough to trivialize a moduli point, for $a$ large enough depending only on $X$ and the rank of $\sE$ (\cite[Prop.~3.2]{ES15}). One applies this to the Frobenius pull-backs of $E_X$. For this one needs that they are semistable as well, and this is the reason for the assumption on $\mu_{\rm max}(\Omega^1_X)$. On the other hand, a Langton type argument guarantees that one finds a crystal $E$ such that $E_X$ is semi-stable. The issue is then to show  vanishing of the Chern classes of $E_X$. It is easy to show that all lattices $E_X$ of a given isocrystal $\sE$ have the same crystalline Chern classes $c^{\rm crys}_n(\sE)$ in $H^{2n}(X/W), n\ge 1$, and thus $c^{\rm crys}_n(\sE)=0$ if $\sE\in {\rm Conv}(X/K)$ (\cite[Prop.~3.1]{ES15}). If $E_X$ is locally free as a coherent sheaf, it is true, but by no means trivial, that $c^{\rm crys}_n(\sE)=0$, where $\sE$ is the isocrystal class of $E$ (see \cite[Section~2/3]{ES15b}). However, we do not know whether or not any isocrystal $\sE$ admits a lattice $E$ which is locally free. This explains the restriction on the type of isocrystals considered in the theorem.

\end{document}